\documentclass[10pt,letterpaper,twocolumn]{article}
\usepackage[utf8]{inputenc}
\usepackage{amsmath}
\usepackage{amsfonts}
\usepackage{amssymb}
\usepackage{xcolor}

\newcommand\R{\mathbb{R}}

\newcommand\T{\mathbb{T}}

\newtheorem{theorem}{Theorem}

\newtheorem{conjecture}[theorem]{Conjecture}

\newcommand\blfootnote[1]{%
  \begingroup
  \renewcommand\thefootnote{}\footnote{#1}%
  \addtocounter{footnote}{-1}%
  \endgroup
}

\title{On turbulence and geometry: from Nash to Onsager}
\author{Camillo De Lellis, L\'aszl\'o Sz\'ekelyhidi Jr.}

\begin{document}

\maketitle

\blfootnote{Communicated by Daniela De Silva\\
Camillo De Lellis is professor of mathematics at the Institute for Avdanced Study and at the University of Z\"urich. His email address is camillo.delellis@math.ias.edu\\
L\'aszl\'o Sz\'ekelyhidi Jr. is professor of mathematics at the University of Leipzig. His email address is laszlo.szekelyhidi@math.uni-leipzig.de}

\subsection*{Turbulence: a challenge for mathematicians}

There is a huge literature on turbulent incompressible flows in applied mathematics, physics and engineering. The outcome of such tremendous effort has been the derivation of several theories, which often allow quite accurate predictions of many phenomena. There is also a quite broad consensus on which fundamental partial differential equations (in short PDEs) describe with sufficient accuracy incompressible fluids, namely the Navier Stokes and the Euler equations. Therefore, in a certain statistical or averaged sense, the predictions of the ``theory of turbulence'' should ultimately translate in mathematically verifiable claims about the behavior of solutions to the latter well-known PDEs. Indeed the working mathematician, even if not immediately confronted with clear cut mathematical statements about turbulence (at least as we understand them in pure mathematics), will be nonetheless able to derive some precise mathematical  consequences of the discussions right at the start of most textbooks about turbulent flows. 

Yet, it seems very hard to prove any of these statements rigorously. Many valuable works in pure mathematics have shown the validity of some of the mechanisms identified by applied mathematicians, physicists and engineers. However these results are mostly confined either to models or to situations obeying some important (and as of yet unverifiable) apriori assumptions. In this note we want to report on a series of recent works which culminated in the complete verification of a famous statement in the theory of turbulence, not confined to some model or constrained by some special a priori assumption. These results have uncovered a surprising and interesting connection with a classical area of differential geometry and have also been used in other contexts. 

\subsection*{K41 theory and the Onsager conjecture}\label{s:K41}

Consider the incompressible Navier-Stokes equations
\begin{equation}\label{e:NavierStokes}
\left\{
\begin{array}{l}
\partial_t v+{\rm div} (v\otimes v) +\nabla p -\nu \Delta v=0\\ \\
\textrm{div }v =0\, ,
\end{array}\right.
\end{equation}
describing the motion of an incompressible homogeneous viscous fluid, where we set the density to be $1$ for simplicity. The pair $(v,p)$ consists, respectively, of a vector and a scalar function: $v (x,t)$ is the velocity of the fluid particle which occupies the point $x$ at time $t$, and $p$ is the hydrodynamic pressure.  Whilst acknowledging that most of the interesting hydrodynamic phenomena in reality happen in the presence of boundaries, in order to isolate a manageable mathematical situation we will assume in the rest of this note that spatial domain is $\mathbb T^3=(\R/2\pi \mathbb Z)^3$, the $3$-dimensional flat torus with side-length $2\pi$. In other words, we consider the problem \eqref{e:NavierStokes} with periodic boundary conditions in the box $[0,2\pi]^3$. 

The coefficient $\nu>0$ is called the kinematic viscosity of the fluid. If the characteristic scale (which we have set to be unity) and the characteristic velocity of the flow are both fixed, then $\nu$ is inversely proportional to the {\em Reynolds number $\rm Re$}, a quantity without physical dimensions, which is characteristic for the turbulent nature of the flow under consideration. A good intuition for the nature of ${\rm Re}$ is given by the following simple model. Replace the left hand side of the first equation in \eqref{e:NavierStokes} with an external force $f$ acting on low wavenumbers, for instance $f(x)=A\sin(x_1)e_2$, where $A>0$ is a parameter and $e_2$ is the unit vector in the $x_2$-direction. In this case the system has a simple stationary solution, given by $\bar{v}(x)=\frac{A}{\nu}\sin(x_1)e_2$. Numerical experiments show however that, as we increase the parameter $A$, the stationary solution eventually becomes unstable, and going through a series of bifurcations the observed flow becomes more and more complex, at some stage becoming chaotic - in other words turbulent.  
Note that in this situation the characteristic velocity is $\frac{A}{\nu}$ and the Reynolds number is given by $Re=\frac{A}{\nu^2}$, so that increasing $A$ has the same effect as decreasing $\nu$. This can be seen by a simple rescaling of time: given a solution $v=v(x,t)$ of \eqref{e:NavierStokes} with the external force $f$ above, set $u(x,t)=\nu^{-1} v(x,\nu^{-1} t)$. Then $u$ is also a solution of the Navier-Stokes system on $\T^3$, with kinematic viscosity equal to $1$ and external force $\nu^{-2} f(x)=Re\,\sin(x_1)e_2$. The simple stationary solution becomes $\bar{u}(x)=\nu^{-1}\bar{v}(x)=Re\,\sin(x_1)e_2$, hence $Re$ is the only parameter remaining. 

Going back to the Navier-Stokes equations with no driving external force, elementary calculations show that for smooth solutions of \eqref{e:NavierStokes} the total kinetic energy
\[
E (t) := \frac{1}{2} \int |v (x,t)|^2\, dx\, 
\]
satisfies the energy balance law 
$$
\frac{d}{dt}E(t)=-\nu\int |\nabla v|^2\,dx.
$$
Thus, at least formally, one would expect that as $\nu\to 0$, the energy dissipation rate vanishes. This naive expectation is contradicted by observation, both physical and numerical: the dissipation rate remains finite and positive. This effect is called {\em anomalous dissipation} in the literature. 
Kolmogorov in the early 1940s pioneered the statistical theory of turbulent motions, assuming that generic flows can be seen as realizations of random fields. Kolmogorov's theory postulates (cf. \cite[Chapter 5]{FrischBook}) that the energy dissipation is strictly positive and independent of the viscosity $\nu$ when the latter goes to $0$ - in agreement with observation. The key insight is that anomalous dissipation arises from the fact that no matter how small the kinetic viscosity $\nu$ is, there is a steady flow or \emph{cascade} of energy from low to high frequencies, leading to large $|\nabla v|^2$. In the words of J.~von Neumann \cite{VonNeumann:35DN87M4}, the decisive trait is that \emph{turbulence is not a matter of ergodic distribution of a fixed amount of energy, but the transport of a fixed flow of energy from sources in the low frequencies to sinks in the high frequencies in the Fourier- transform space}
Assuming in addition local homogeneity and isotropy, Kolmogorov derived his famous $k^{-5/3}$ law, expressing the mean distribution of kinetic energy density across an intermediate range of spatial frequencies $k_0<k<k_{\nu}$. 

When $\nu\downarrow 0$, \eqref{e:NavierStokes} becomes formally the incompressible Euler equations
\begin{equation}\label{e:Euler}
\left\{
\begin{array}{l}
\partial_t v + {\rm div} (v\otimes v) + \nabla p = 0,\\ \\
{\rm div}\, v = 0\, . \\
\end{array}\right.
\end{equation}
The energy balance law implies that smooth solutions of \eqref{e:Euler} preserve the total kinetic energy.

Onsager suggested in his famous note \cite{Onsager:1949kz} the possibility of anomalous dissipation for {\it weak solutions} of the Euler equations as a consequence of the energy cascade. Indeed, at least formally, as $\nu\downarrow 0$, the inertial range of frequencies extends to infinity (namely $k_{\nu}\to \infty$), hence Kolmogorov's $k^{-5/3}$ law amounts to a certain regularity statement, when interpreted for single velocity fields rather than ensemble averages. This is exactly what Onsager proposed in 1949\footnote{Towards the end of his note Onsager writes:

{\em ... It is of interest to note that in principle, turbulent dissipation as described could take place just as readily without the final assistance of viscosity {\rm [A/N: in the previous pages Onsager discusses the energy spectrum of turbulent solutions of \eqref{e:NavierStokes} on the threedimensional torus after rewriting the equations as an infinite-dimensional system of ODEs for the Fourier coefficients; the ``absence of viscosity'' refers to setting $\nu =0$ in \eqref{e:NavierStokes}]}. In the absence of viscosity, the standard proof of the energy conservation does not apply, because the velocity field does not remain differentiable! In fact it is possible to show that the velocity field in such ``ideal'' turbulence cannot obey any Lipschitz condition of the form
\[
|\vec{v} (\vec{r}' + \vec{r}) - \vec{v} (\vec{r}')| < ({\rm const.}) r^n
\]
for any $n$ greater than $\frac{1}{3}$; otherwise the energy is conserved.}}. It is important to emphasize that the theory of Kolmogorov is a {\it statistical theory}, dealing with random fields whose distribution laws need to satisfy several postulates, but even the mere existence of such random fields is not at all obvious in rigorous mathematical terms. In contrast, the suggestion of Onsager turned the problem into a ``pure PDE'' question that could be studied directly, and, after nearly 70 years, we can finally state the theorem:

\begin{theorem}\label{c:Onsager}
Let $(v,p)$ be a weak solution of \eqref{e:Euler}  with
\begin{equation}\label{e:uniformholder}
|v(x,t)-v(y,t)|\leq C|x-y|^{\theta} \qquad \forall x,y,t
\end{equation}
(where $C$ is a constant independent of $x,y,t$). 
\begin{itemize}
\item[(a)] If $\theta> \frac{1}{3}$, then $E(t)$ is necessarily constant; 
\item[(b)] For $\theta< \frac{1}{3}$ there are solutions for which $E (t)$ is strictly decreasing. 
\end{itemize}
\end{theorem}

In the historical context it is quite remarkable that Onsager was able to formulate a mathematically very precise statement. He gave, in particular, a rigorous definition of ``weak solutions'' expanding \eqref{e:Euler} into an infinite system of ODEs for its Fourier coefficients\footnote{Onsager's continuation of the aforementioned paragraph reads: 

{\em Of course, under the circumstances, the ordinary formulation of the laws of motion in terms of differential equations becomes inadequate and must be replaced by a more general description; for example the formulation (15) {\rm [A/N: the reference is to formula (15) in Onsager's paper]} in terms of Fourier series will do.}}, a point of view which coincides with the concept of ``distributional solution'' of the modern PDE literature. Nevertheless, it was only in the early 1990s that mathematicians took note of this statement as a mathematical conjecture, mostly as a result of Greg Eyink's efforts in providing a modern account of Onsager's unpublished work on the topic \cite{Eyink:2006jn} and in popularizing the subject in the mathematical fluid dynamics community \cite{Eyink:1994vd}.  

Part (a) of Theorem \ref{c:Onsager} was proved in \cite{Constantin:1994vn} using a regularization procedure and a clever and powerful, yet elementary, commutator estimate.

Part (b) of Theorem \ref{c:Onsager} took another 25 years, with a series of partial results and gradual improvements \cite{Scheffer93,Shnirelman:1997uz,Shnirelman:2000vu,DeLellis:2009jh,DeLellis:2008vc,DeLellis:2013im,DeLellis:2012tz,Isett:2013ux,Buckmaster:wb,Buckmaster:2013vv,Buckmaster:2014ty,SzekelyhidiJr:2016tp}, finally culminating in \cite{Isett:2016to} and the subsequent improvement \cite{Buckmaster:2017uza}.

Well-formulated mathematical conjectures are not just about solving a problem. We associate with them the hope of learning something deeper about the context in which the problem was formulated, possibly revealing unexpected connections between different parts of mathematics. The story of Onsager's conjecture fits very well in this ideal: the proof of Theorem \ref{c:Onsager} is quite possibly more interesting than the actual statement (which, in some sense, was already ``clear'' to physicists). It reveals a new mechanism by which energy cascades may appear in solutions of nonlinear PDEs. This mechanism has already led to further breakthroughs, see below, and has uncovered an entirely unexpected connection between the mathematically elusive concept of turbulence and the ``New Land" (cf. \cite{Gromov:2015tua}) discovered by John Nash in differential geometry.

\subsection*{The Scheffer-Shnirelman paradox}

In \cite{Scheffer93} V.~Scheffer constructed a non-trivial weak solution to the Euler equations in $\R^2$ with compact support in space and time. Subsequently A.~Shnirelman in \cite{Shnirelman:1997uz} gave an entirely different proof in $\T^2$. Such a result is hard to interpret physically, as it would correspond to a perfect incompressible fluid which can start and stop moving by itself, without the action of external forces. As such, for a long time Scheffer's Theorem remained some sort of ``paradox'' in the PDE literature \cite{Villani:2008vp}, cited mostly as a warning example of unphysical behaviour if the notion of solution is too weak, with emphasis more on the \emph{non-uniqueness} aspect rather than the \emph{violation of energy conservation}. Indeed, since the weak solutions constructed here are merely square-integrable in space-time, there is no control on the regularity of the total kinetic energy $E(t)$ (which is not even known to be bounded for every $t$), therefore in principle $E$ need not be monotone on any time-interval. The first result in connection with Theorem \ref{c:Onsager}(b) was obtained by Shnirelman in \cite{Shnirelman:2000vu} - he showed the existence of a solution on $\mathbb T^3$ with finite and strictly monotone decreasing energy on some time interval. Nevertheless, this solution is not continuous, and, since it was obtained using a generalized flow model with sticky particles, it seemed to have no connection to the energy cascade picture postulated by Kolmogorov and Onsager. 

\subsection*{Differential inclusions, relaxation and Reynolds stress}

Much less known than his works on the Euler and Navier-Stokes equations, the unpublished PhD thesis of Scheffer \cite{Scheffer_PhD} contains some remarkable precursors of fundamental results in the vectorial calculus of variations.  One of them is akin to Evans' $\varepsilon$-regularity theorem \cite{Evans:1986up} for minimizers of elliptic energies (the precise assumption is ``uniform quasiconvexity'') and the other one is a precursor of a striking singularity theorem by M\"uller and \v{S}verak \cite{Muller:2003wy}. The latter paper (which as we will see below has been an important source of inspiration for at least another reason) proves the existence of a rather badly behaved Lipschitz critical point to elliptic energies which fall under the assumptions of Evans' theorem: such critical point is nowhere differentiable, in stark contrast with minimizers, which according to Evans must be smooth on a dense open subset. 

The work \cite{Muller:2003wy} belongs, unlike the papers \cite{Scheffer93,Shnirelman:1997uz,Shnirelman:2000vu}, to an established mathematical tradition, in the sense that during the 1990s and early 2000 several authors used analogous ideas to produce rather striking examples of irregular solutions to various systems of partial differential equations, all falling in the class of ``differential inclusions'' (the most notable examples are the works \cite{Dacorogna:1997wk,Sychev:2001wo,Kirchheim:2003wp}). In fact an important functional analytic aspect of these ideas was pioneered in the context of ordinary differential equations in the work \cite{Cellina}, see also \cite{BressanFlores,Cellina:2005tx}. The basic principle underlying these works can be explained on the simplest of differential inclusions: consider the sets
\begin{equation}\label{e:1D-DI}
\begin{split}
X&=\bigl\{z:[0,1]\to\R:\,|z|\leq 1\textrm{ a.e.}\bigr\}\\
S&=\bigl\{z:[0,1]\to\R:\,|z|=1\textrm{ a.e.}\bigr\}\,.
\end{split}
\end{equation}
Then $S$ is dense in $X$ with respect to the weak* topology in $L^{\infty}$ (in fact Baire-residual), see \cite{Cellina,BressanFlores}. This statement easily follows from Mazur's lemma and emphasizes the relationship between weak convergence and convex hulls. 

The jump from this basic principle to partial differential equations was made possible by the seminal work \cite{Tartar:1977vo} by L.~Tartar (see also \cite{DiPerna:1985vb} for analogous ideas appearing at the same time), which examined the relations between weak convergence and differential constraints and uncovered new phenomena depending on the core PDE structure. As an example, consider the following $n$-dimensional generalization of \eqref{e:1D-DI}:
\begin{equation}\label{e:nD-DI}
\begin{split}
X&=\bigl\{u:B\to\R^n:\,\nabla u\in \textrm{co}K\textrm{ a.e.}\bigr\}\\
S&=\bigl\{u:B\to\R^n:\,\nabla u\in K\textrm{ a.e.}\bigr\}\,,
\end{split}
\end{equation}
where $B\subset\R^n$ is the open unit ball, $K$ is a compact subset of $n\times n$ matrices and $\textrm{co}K$ denotes the convex hull of $K$. 
If $K=O(n)$, the set of linear isometries, then $S$ is dense in $X$ with respect to the uniform topology \cite{Gromov}. However, if $K=SO(n)$, the set of orientation-preserving linear isometries, then $S$ consists only of affine functions \cite{Reshetnyak:1994ts}. The situation for a general $K$ lies somewhere in between these extreme cases, see \cite{Muller:1998wh,Kirchheim:2002wc}.
  
It turns out \cite{DeLellis:2009jh} that the Euler system \eqref{e:Euler}, when interpreted as a differential inclusion similar to \eqref{e:nD-DI}, falls in the same category as the $O(n)$-case above. In order to explain this, note that the density of $S$ in $X$ can be interpreted as a relaxation statement. Consider a sequence of Lipschitz maps $u_k:B\to\R^n$ with $\nabla u(x)\in O(n)$ for a.e. $x\in B$, i.e.
\begin{equation}\label{e:On}
\partial_iu\cdot\partial_ju=\delta_{ij}\textrm{ a.e.},
\end{equation}
which converge uniformly to a limit $u$. Without any additional information on the sequence, we can only say that the limit satisfies $\nabla u(x)\in \textrm{co}O(n)$ a.e., i.e. 
\begin{equation}\label{e:coOn}
\partial_iu\cdot\partial_ju\leq \delta_{ij}\textrm{ a.e.}
\end{equation}
in the sense of quadratic forms. In more geometric language, the uniform limit of a sequence of isometries is a short map. The density of $S$ in $X$ is a kind of converse of this statement, that is, the system \eqref{e:coOn} is the \emph{relaxation} of the system \eqref{e:On}.

The analogous question for the Euler equations is as follows: Consider a sequence of uniformly bounded weak solutions $v_k$ to \eqref{e:Euler}, which converge weakly in $L^2$ to some limit function $\bar v$. Then $\bar v$ satisfies the Euler equations with an error term:
\begin{equation}\label{e:Reynolds}
\left\{\begin{array}{l}   
\partial_t\overline{v}+{\rm div}(\overline{v}\otimes\overline{v}+ R)+\nabla\overline{p} =0\\ \\
{\rm div}\overline{v} =0,
\end{array}\right.
\end{equation}
where, denoting by $\overline{v\otimes v}$ the weak limit of $v_k\otimes v_k$, 
\begin{equation}\label{e:R}
R=\overline{v\otimes v}-\overline{v}\otimes\overline{v}\, .
\end{equation}
Without any additional information on the sequence $v_k$, we have no reason to expect that the symmetric tensor $R$ vanishes, the only information we can gather is that $R(x,t)\geq 0$ a.e. in the sense of quadratic forms (the passage from $R=0$ to $R\geq 0$ is the precise analogue to the passage from $K$ to $\textrm{co}K$ in \eqref{e:nD-DI}). 
The tensor $R$ is, in a slightly different guise, a very well known object in the theory of turbulence, called Reynolds stress, which arises by a formal averaging of the equations \eqref{e:NavierStokes} or \eqref{e:Euler}. Thinking of a turbulent velocity field as the sum $v=\bar{v}+w$ of a mean flow $\bar{v}$ and a (random) fluctuation, the induced stress by the fluctuation on the mean flow is given by $R=\overline{w\otimes w}$. This is the same formula as \eqref{e:R}. Indeed, the process of taking weak limits is somewhat akin to averaging high-frequency oscillations, and weak convergence in place of averaging random fluctuations has been proposed by P. Lax \cite{Lax:1991ve} as a deterministic approach to turbulence.

The work in \cite{DeLellis:2009jh,DeLellis:2008vc} established the relationship between the system \eqref{e:Reynolds} and \eqref{e:Euler}. More precisely, 
the main result in \cite{DeLellis:2008vc} states that (modulo some technical assumptions) any solution $\bar{v}$ of \eqref{e:Reynolds}  can be weakly approximated by weak (bounded, but in general discontinuous) solutions $v$ of \eqref{e:Euler}, so that, in this sense, the system \eqref{e:Reynolds} is the relaxation of the Euler system. In this way we were able to place the existence theorems of Scheffer and Shnirelman in a very general and flexible context, applicable not just to the incompressible Euler equations, but to several other PDEs, see e.g. \cite{Cordoba:2011wp,Shvydkoy:2011ta,Chiodaroli:2014cw} (note in contrast, that Shnirelman's proofs heavily relied on generalized flows and thus on Arnold's geodesic flow formulation of the incompressible Euler equations). Moreover, our work provided a surprising {\it fil rouge} between Scheffer's PhD thesis and his work \cite{Scheffer93} on the incompressible Euler equations. 

Note however that, in terms of regularity and the energy cascade picture, these results were no closer to Theorem \ref{c:Onsager}(b) than \cite{Scheffer93,Shnirelman:1997uz,Shnirelman:2000vu}.

\subsection*{Convex integration and the Nash-Kuiper paradox}

As mentioned above, there is a second reason why the paper \cite{Muller:2003wy} by S.~M\"uller and V.~\v{S}verak was a major source of inspiration for us. The authors in \cite{Muller:2003wy} pointed out for the first time an important relation between the results in the theory of differential inclusions and Gromov's $h$-principle in geometry. In particular the method of convex integration, introduced by M.~Gromov \cite{Gromov} and extended
in \cite{Muller:2003wy} to Lipschitz mappings, provides a very
powerful tool to construct solutions to nonlinear PDEs. 

The origin of Gromov's convex integration lies in the famous Nash-Kuiper theorem on isometric embeddings of Riemannian manifolds.
Let $\Sigma$ be a smooth compact manifold 
of dimension $n\geq 2$, equipped with a Riemannian metric $g$.  A map $u : \Sigma \to \mathbb R^N$ is {\em isometric} if it preserves the length of curves, i.e. if
\begin{equation}\label{e:true_isometry}
\ell_g (\gamma) = \ell_e (u\circ \gamma) \qquad \mbox{for any $C^1$ curve $\gamma\subset \Sigma$,}
\end{equation}
where $\ell_g (\gamma)$ denotes the length of $\gamma$ with respect to the metric $g$:
\begin{equation}\label{e:Riem_length}
\ell_g (\gamma) = \int \sqrt{ g (\gamma (t)) [\dot\gamma (t), \dot\gamma (t)]}\, dt\, .
\end{equation}
If $u\in C^1(\Sigma;\mathbb R^N)$ this means, using the language of Riemannian geometry, that the pull back of the Euclidean metric $u^\sharp e$
agrees with $g$. In local coordinates such relation is the following system of partial differential equations
\begin{equation}\label{e:equations}
\partial_iu\cdot\partial_ju=g_{ij}\, ,
\end{equation}
which can be seen as an (inhomogeneous) differential inclusion - compare with \eqref{e:On}. 

The existence of isometric immersions (and embeddings) of Riemannian manifolds into some Euclidean space is a classical problem, explicitly formulated for the first time by Schl\"afli, see \cite{Schlaefli}, who conjectured that the system is solvable {\em locally}
if the dimension $N$ of the target is at least $\frac{n(n+1)}{2}$, matching the number of equations in \eqref{e:equations}. An isometric immersion in co-dimension $1$ would seem a rare bird, since, with the exception of the case $n=2$, the system \eqref{e:equations} would be heavily overdetermined. Yet, in \cite{Nash:1954vt} J.~Nash astonished the geometry world by proving essentially that the only obstruction to the existence of solutions of \eqref{e:equations} is topological, at least in the class of $C^1$ maps. 

In order to state Nash's Theorem let us recall that an immersion $u: \Sigma \to \mathbb R^N$ is called short if it ``shrinks'' the length of curves. For $C^1$ immersions and in local coordinates such condition is equivalent to the inequality 
\begin{equation}\label{e:sub_for_nash}
\partial_iu\cdot\partial_ju \leq g_{ij}
\end{equation}
in the sense of quadratic forms - compare with \eqref{e:coOn}.

\begin{theorem}\label{t:main_C1_1}
Let $(\Sigma, g)$ be a smooth closed $n$-dimensional Riemannian manifold. Any $C^\infty$ short immersion $u: \Sigma \to \mathbb R^N$ with $N\geq n+1$ can be uniformly approximated by $C^1$ isometric immersions. If $u$ is, in addition, an embedding, then it can be approximated by $C^1$ isometric embeddings. 
\end{theorem}

Nash proved Theorem \ref{t:main_C1_1} for $N\geq n+2$ and suggested that his strategy could be suitably modified to work in the case $N=n+1$; the details were then given in two subsequent works by N.~Kuiper \cite{Kuiper:1955th}. 

Theorem \ref{t:main_C1_1} can be seen as the $C^1$ analogue of the relaxation statement in \eqref{e:nD-DI}. However, whilst it is rather easy to imagine Lipschitz isometric maps arising as ``foldings'' of the manifold $\Sigma$, the $C^1$ case came as a complete surprise, in particular because the Nash-Kuiper theorem cannot hold for $C^2$ maps. For instance, a $C^2$ isometric immersion of a closed positively curved sphere in the three-dimensional space is necessarily convex, and in fact the shape is determined up to rigid motions by a classical result of Cohn-Vossen and Herglotz \cite{CohnVossen:1930iy,Herglotz:1943je}. In the 1950s Yu.~Borisov showed that such result can be extended to $C^{1,\frac{2}{3}+\varepsilon}$ isometries, cf. \cite{Borisov:1959tt}. In fact, since for a $C^1$ surface the Gauss map is continuous, with a well-defined Brouwer degree, there was some hope that the rigidity statement of Cohn-Vossen and Herglotz can be extended to to the $C^1$ case. It was this hope that was shattered by Nash's result.   

Subsequently, Borisov announced that the Nash-Kuiper theorem can be extended to $C^{1,\alpha}$ isometries provided $\alpha$ is sufficiently small \cite{Borisov:1965wf}. While a detailed proof of these announcements only appeared in one special case in \cite{Borisov:2004wo}, in the joint work \cite{Conti:2012wl} with Sergio Conti we revisited these extensions and provided a unified framework for all the results announced in \cite{Borisov:1965wf}. In that paper we also noticed that the geometric considerations leading to Borisov's rigidity statement \cite{Borisov:1959tt} can be substituted by a short PDE argument which relies on the same commutator estimate as in Constantin, E and Titi's proof of Theorem \ref{c:Onsager}(a). 

Thus, a striking analogy between isometric immersions and solutions of the Euler equations arose: at sufficiently high regularity ($C^2$ for isometries, $C^1$ for Euler) solutions are well-behaved and with appropriate side-conditions uniquely determined, at sufficiently low regularity (Lipschitz/$C^1$ for isometries, bounded/$C^0$  for Euler) solutions with completely different behavior appear. It was thus natural to try to adapt the ideas of Nash in \cite{Nash:1954vt} to the Euler equations. It is important to emphasize, however, that this analogy concerns more the \emph{non-uniqueness} aspect of weak solutions of Euler, i.e. the unphysical behavior already observed in connection with the Scheffer-Shnirelman paradox. On the other hand, in light of Gromov's h-principle, one should perhaps view this aspect more as an expression of \emph{flexibility} rather than non-uniqueness - the violation of energy conservation is one aspect of this flexibility. Indeed, while part (a) of Theorem \ref{c:Onsager} shows that above the H\"older exponent $\frac{1}{3}$ solutions cannot be too flexible, it is not at all clear what to expect about their uniqueness: one guess might be that the threshold for uniqueness is a small improvement of $C^1$ (in the Osgood sense).

\subsection*{Dissipative continuous solutions}

Inspired by Nash's proof in \cite{Nash:1954vt} we devised in \cite{DeLellis:2013im} a ``convex integration'' scheme leading to continuous dissipative solutions of \eqref{e:Euler}. Subsequently we showed in \cite{DeLellis:2012tz} that these solutions satisfy Theorem \ref{c:Onsager}(b) with exponent $\theta<\frac{1}{10}$. The construction is based on an 
iteration, where at each step we add a highly oscillatory correction in order to decrease the defect to being a solution. More precisely, we construct inductively a sequence of solutions $(v_q,p_q,R_q)$, $q=1,2,\dots$ to
\begin{equation}\label{e:Reynolds_2}
\left\{\begin{array}{l}
\partial_t v_q+{\rm div}(v_q\otimes v_q)+ \nabla p_q = - {\rm div}\, R_q\,,\\ \\
{\rm div} v_q =0\,, 
\end{array}\right.
\end{equation}
such that $v_q\to v$ and $R_q\to 0$ uniformly. Observe that \eqref{e:Reynolds_2} is the same system as \eqref{e:Reynolds}. Accordingly, $v_{q+1} = v_q+w_{q+1}$ where we think of $v_q$ as the ``mean flow'' on length-scales $\geq \lambda_q^{-1}$ and $w_{q+1}$ is the ``fluctuation'' on this scale. Thus, up to lower order corrections $w_{q+1}$ should have the form
\begin{equation}\label{e:ansatz_1}
w_{q+1}(x,t)=W\Bigl(v_q (x,t),  R_q (x,t),\lambda_{q+1} x,\lambda_{q+1} t\Bigr)\, ,
\end{equation}
where $W (v, R, \xi, \tau)$ is some ``master function'' and $\lambda_{q+1}$ a parameter which increases at least exponentially fast at each step. In comparison, in the proof of Nash \cite{Nash:1954vt} these ``fluctuations'' are spirals (and in \cite{Kuiper:1955th} corrugations) aimed at increasing the metric -- thereby reducing the metric error -- in a single coordinate direction. 

The basic idea for reducing the error with such an \emph{Ansatz} is the following: assuming that $v_q$ is already the correct solution up to spatial frequencies of order $\lambda_q$, and $w_{q+1}$ is supported on spatial frequencies of order $\lambda_{q+1}$, it is easy to see that the only possibility for $w_{q+1}$ to correct the error $R_q$ is via the high-high to low interaction in the product $w_{q+1}\otimes w_{q+1}$. In other words, the master function $W$ should satisfy the properties that $\xi\mapsto W(v,R,\xi,\tau)$ is $2\pi$-periodic with average 
\begin{equation*}
\langle W\rangle := \frac{1}{(2\pi)^3}\int_{\T^3}W(v,R,\xi,\tau)\,d\xi=0;
\end{equation*}
and the average stress is given by $R$, i.e. 
\begin{equation*}
\langle W\otimes W\rangle = R\, ;
\end{equation*}
Note how these requirements are consistent with \eqref{e:R}. 
Now assume that $\xi\to W(v,R,\xi,\tau)$ is a stationary solution of Euler for any $v,R,\tau$. Substituting this \emph{Ansatz} for $v_{q+1}$ into \eqref{e:Reynolds_2} then yields as the main quadratic interaction term ${\rm div} (w_{q+1}\otimes w_{q+1}-R_q)+\nabla p_{q+1}$, leading to a new Reynolds stress $R_{q+1}$ with Fourier support on frequencies of order $\lambda_{q+1}$. In this way we can push the error to high frequencies by successively ``undoing'' the averaging process leading to Reynolds stresses in \eqref{e:Reynolds}-\eqref{e:R}. 

As explained in \cite{SzekelyhidiJr:2016wc,DeLellis:2017dt}, starting from the \emph{Ansatz} above it is possible to write down a family of conditions that $W$ would have to satisfy, ideally, so to give a ``clean'' convex integration iteration leading to a proof of Theorem \ref{c:Onsager}(b). Although this family of conditions is somewhat naive and unfortunately no such $W$ exists (indeed, the scaling of time in \eqref{e:ansatz_1} is clearly ``wrong''), approximations based on a special family of stationary solutions of Euler called Beltrami flows can be used. This lead to the results in \cite{DeLellis:2013im,DeLellis:2012tz}.

\subsection*{Climbing the Onsager ladder}

The reason why the construction in \cite{DeLellis:2013im,DeLellis:2012tz} was only able to produce weak solutions to the Euler system as in Theorem \ref{c:Onsager}(b) with H\"older exponent $\theta<\frac{1}{10}$ is the rather poor control of the linear (i.e.~transport) interaction term between ``mean flow'' $v_q$ and ``fluctuation'' $w_{q+1}$. Indeed, whilst it is quite clear from dimensional considerations that the scaling $(\lambda_{q+1}x,\lambda_{q+1}t)$ in \eqref{e:ansatz_1} is unnatural, several modifications were introduced later in the precise implementation of the basic iteration scheme described above. P.~Isett introduced in his PhD thesis \cite{Isett:2013ux} the correct space-time scaling, which eventually lead to an improvement of Theorem \ref{c:Onsager}(b) with H\"older exponent $\theta<\frac{1}{5}$, cf. \cite{Buckmaster:2014ty}. Moreover, following an idea from the PhD thesis of T.~Buckmaster to introduce temporal intermittency \cite{Buckmaster:2013vv,Buckmaster:wb}, in the joint work \cite{Buckmaster:2014th} we reached for the first time the threshold $\theta<\frac{1}{3}$, although not in the desired scale of spaces: more precisely, for any $\theta<\frac{1}{3}$ we show the existence of nontrivial continuous solutions with compact temporal support (i.e.~in the spirit of the Scheffer-Shnirelman paradox) which satisfy the condition
\[
|v (x,t) - v(y,t)|\leq C (t) |x-y|^{\theta} \qquad \mbox{for every $x,y,t$}\, ,
\]
where $C$ is an $L^1$ function of time.

The correct scale of spaces was finally achieved in \cite{Isett:2016to}, where P.~Isett proved the existence of compactly supported nontrivial solutions in $C^{\theta}$ for every $\theta<\frac{1}{3}$. The proof in \cite{Isett:2016to} contains two new ideas. Firstly, in \cite{SzekelyhidiJr:2016tp} S.~Daneri and the second author introduced a new class of ``master functions'' $W$ called Mikado flows, which are more stable under convection by a large-scale mean flow. However, it is not possible to use such $W$ directly in the scheme of \cite{DeLellis:2013im} (and its subsequent refinements), not even to produce continuous solutions. The second key idea in \cite{Isett:2016to} is a gluing technique, which combines the convex integration technique with the free (unforced) Euler dynamics in an alternating fashion. A shortcoming of \cite{Isett:2016to} is the insufficient control of the kinetic energy $E(t)$. The final proof of Theorem \ref{c:Onsager}(b) was then given in \cite{Buckmaster:2017uza}.

\subsection*{Further developments and open problems}

Several interesting challenges remain in the area. One is to produce a sequence of Leray-Hopf solutions to the Navier-Stokes equations with vanishing viscosity which exhibit anomalous dissipation. This would be the case if, for instance, one were able to show that at least one of the solutions of Euler constructed by the convex integration method is a limit of Leray-Hopf solutions of the Navier-Stokes equations with vanishing viscosity. The latter problem might be linked to constructing dissipative solutions which are Onsager-critical, although it is not clear which scale of spaces is most natural.

A second important challenge in connection with the Kolmogorov-Onsager theory is to introduce intermittency and in this way to be able to recover the measured deviations from K41 self-similarity and possible multifractality. In their remarkable recent work \cite{Buckmaster:2017wfa} Buckmaster and Vicol have constructed a convex integration scheme which produces {\em weak} solutions of the Navier-Stokes equations, by introducing intermittency. In particular they were able to show that:
\begin{itemize}
\item For such weak solutions the Cauchy problem to the Navier-Stokes equations is ill-posed;
\item Any ``convex-integration solution'' of Euler can be approximated strongly in $L^2$ by weak solutions of Navier-Stokes.
\end{itemize}
Furthermore, in the joint work \cite{Buckmaster:2018wj} with M.~Colombo they could produce such solutions with the additional property that they are smooth except for a small closed set of times (of Hausdorff dimension strictly smaller than $1$). While the above results are striking, they seem for the moment far from producing solutions that belong to the energy space. The latter can be instead reached when the Laplacian in the viscosity is replaced by a fractional Laplacian with sufficiently low exponent
(cf. \cite{CDR,DeRosa}; the latter works follow more closely the constructions in \cite{Buckmaster:2014ty} and \cite{Buckmaster:2017uza}).
In a similar line of research, the second author in joint work with S.~Modena \cite{Modena:2017tp,Modena:2018tl} constructed convex integration solutions to the transport and continuity equations with Sobolev vector fields which are not renormalized, thus showing optimality of the integrability assumptions in the DiPerna-Lions theory.  

A further important issue is related to the closure problem and dynamical behavior of the Reynolds stress. While the results in \cite{DeLellis:2008vc,SzekelyhidiJr:2016tp,Buckmaster:2017uza} imply that in the absence of boundaries there is no constraint on the evolution of the Reynolds stress tensor, a series of results on initial value problems associated to classical hydrodynamic instabilities \cite{Szekelyhidi:2011bj,SzekelyhidiJr:2012wo,Castro:2016ut,Forster:2017uw} point at the possibility that in the presence of initial conditions or boundaries this is not the case anymore. The constrained evolution of the Reynolds stress in the general situation remains to be explored.

Regarding the analogy between Euler and isometric immersions, note that Theorem \ref{c:Onsager} provides a sharp threshold between two competing situations. A sort of analog of the Onsager's conjecture has been recently proved for the isometric embedding problem, where the threshold exponent turns out to be $\frac{1}{2}$. In the recent work \cite{DeLellis:2018ub} D.~Inauen and the first author have shown that, while isometric embeddings of class $C^{1,\frac{1}{2}+\varepsilon}$ satisfy a suitable generalization of a classical theorem in differential geometry (namely the Levi-Civita connection of the manifold coincides with the connection induced on the immersion by the ambient Euclidean space), for every $\varepsilon>0$ it is possible to construct $C^{1,\frac{1}{2}-\varepsilon}$ isometric embeddings which violate it.
Even though the latter result shows the criticality of the exponent $\frac{1}{2}$, the most compelling conjecture in the area is still unsolved:

\begin{conjecture}
Consider a positively curved sphere $(\mathbb S^2, g)$ and isometric immersions
$v\in C^{1,\alpha} (\mathbb S^2, \mathbb R^3)$.
\begin{itemize}
\item[(a)] If $\alpha > \frac{1}{2}$, then $v (\mathbb S^2)$ is convex and it is unique up to ambient isometries.
\item[(b)] If $\alpha < \frac{1}{2}$ any short immersion (resp. embedding) $u$ can be uniformly approximated by a sequence of isometric immersions (resp. embedding) $v_k\in C^{1,\alpha} (\mathbb S^2)$.
\end{itemize}
\end{conjecture}

As already mentioned, part (a) of the Conjecture is known to hold for $\alpha >\frac{2}{3}$. Concerning part (b), the joint work of the authors with Inauen \cite{DeLellis:2015wm} prove the statement for positively curved disks when $\alpha < \frac{1}{5}$, while forthcoming work of the second author with W.~Cao settles the case of general surfaces for $\alpha <\frac{1}{5}$.

\bibliographystyle{alpha}


\end{document}